\documentclass[journal]{IEEEtran}
\usepackage{amsmath,amsfonts,amssymb,mathrsfs}
\usepackage{algorithmic}
\usepackage{algorithm}
\usepackage{array}
\usepackage{textcomp}
\usepackage{stfloats}
\usepackage{url}
\usepackage{verbatim}
\usepackage{graphicx}
\usepackage{cite}
\usepackage{booktabs}
\usepackage{caption}
\usepackage{subfigure}
\usepackage{nomencl}
\usepackage{multirow}
\usepackage{threeparttable}
\usepackage{tabularx}
\usepackage{color}
\usepackage{courier}
\allowdisplaybreaks
\makeatletter
\newcommand{\removelatexerror}{\let\@latex@error\@gobble}
\makeatother
\removelatexerror
\renewcommand{\algorithmicrequire}{\textbf{Input:}}
\renewcommand{\algorithmicensure}{\textbf{Output:}}

\begin{document}

\title{Multiple Joint Chance Constraints Approximation for Uncertainty Modeling in Dispatch Problems}

\author{Yilin~Wen,~\IEEEmembership{Graduate Student Member,~IEEE,}
        Yi~Guo,~\IEEEmembership{Member,~IEEE,}
        Zechun~Hu,~\IEEEmembership{Senior Member,~IEEE,}
        and Gabriela~Hug,~\IEEEmembership{Senior Member,~IEEE}
}

\maketitle

\begin{abstract}
  Uncertainty modeling has become increasingly important in power system decision-making. The widely-used tractable uncertainty modeling method---chance constraints with Conditional Value at Risk (CVaR) approximation, can be over-conservative and even turn an originally feasible problem into an infeasible one. This paper proposes a new approximation method for multiple joint chance constraints (JCCs) to model the uncertainty in dispatch problems, which solves the conservativeness and potential infeasibility concerns of CVaR. The proposed method is also convenient for controlling the risk levels of different JCCs, which is necessary for power system applications since different resources may be affected by varying degrees of uncertainty or have different importance to the system. We then formulate a data-driven distributionally robust chance-constrained programming model for the power system multiperiod dispatch problem and leverage the proposed approximation method to solve it. In the numerical simulations, two small general examples clearly demonstrate the superiority of the proposed method, and the results of the multiperiod dispatch problem on IEEE test cases verify its practicality.

\end{abstract}

\section*{Nomenclature}
\subsection*{Abbreviations:}
\begin{IEEEdescription}[\IEEEusemathlabelsep\IEEEsetlabelwidth{DRCC}]
 \item[DER] Distributed energy resource.
 \item[EV] Electric vehicle.
 \item[TCL] Thermostatically controlled load.
 \item[DESS] Distributed energy storage system.
 \item[HVAC] Heating, ventilating, and air conditioning system.
 \item[ADN ] Active distribution network.
 \item[CCP] Chance-constrained programming.
 \item[ICC] Individual chance constraint.
 \item[JCC] Joint chance constraint.
 \item[CVaR] Conditional value at risk.
 \item[DRCC] Distributionally robust chance constraint.
 \item[DRJCC] Distributionally robust joint chance constraint.
\end{IEEEdescription}

\subsection*{Indices and Sets:}
\begin{IEEEdescription}[\IEEEusemathlabelsep\IEEEsetlabelwidth{$\mathcal D_i, \mathcal G_i, \mathcal R$}]
 \item[$t/T$] Index/Number of time slots.
 \item[$i/\mathcal I$] Index/Set of buses, where $i=0$ represents the slack bus.
 \item[$l/{{\mathcal L}}$] Index/Set of branches.
 \item[$d/\mathcal D$] Index/Set of ADNs.
 \item[$g/\mathcal G$] Index/Set of conventional generators.
 \item[$w/\mathcal W$] Index/Set of renewables.
 \item[$\mathcal D_i, \mathcal G_i, \mathcal W_i$] Set of ADNs, conventional generators and renewables connected to bus $i$.
 \item[$s/S_g$] Index/Number of divided segments for cost function of generator $g$.
 \item[$\Xi$] Support set of uncertain variable $\tilde{\boldsymbol\xi}$.
\end{IEEEdescription}
\subsection*{Parameters:}
\begin{IEEEdescription}[\IEEEusemathlabelsep\IEEEsetlabelwidth{$Ptttttt$}]
 \item[$\Delta T$] Length of each time slot.
 \item[$\psi_{l,i}$] Power transfer distribution factor from bus $i$ to line $l$.
 \item[$P_{l}^{\text{U}}$] Transmission capacity of line $l$.
 \item[$P_g^{\text{max/min}}$] Maximum/Minimum output of generator $g$.
 \item[$P_{g,s}^{\text{U}}$] Upper power bound of segment $s$ of generator $g$.
 \item[$r_g^{\text{dn/up}}$] Ramp down/up ability of generator $g$.
 \item[$C_{g,s}$] Energy cost coefficient of output segment $s$ of generator $g$.
 \item[$C_{g}^{+/-}$] Up-/Down-reserve cost coefficient of generator $g$.
 \item[$P_{i,t}^{\text{fix}}$] Fixed load of bus ${i}$ at time slot $t$.
 \item[$P_{w,t}$] Forecasted power of renewable $w$ at time slot $t$. 
 \item[$C_{d}^{+}$] Up-reserve cost coefficient of ADN $d$.
 \item[$\epsilon_{g},\epsilon_{d},\epsilon_{l}$] Risk budgets of the system operator toward generator $g$, distribution network $d$, and transmission line $l$, $\epsilon_g, \epsilon_d, \epsilon_l \in [0,1)$.
 \item[$\rho_{g},\rho_{d},\rho_{l}$] Radii of the Wasserstein ambiguity sets for the DRJCCs of generator $g$, distribution network $d$, and transmission line $l$, $\epsilon_g, \epsilon_d, \epsilon_l\ge 0$.
\end{IEEEdescription}
\subsection*{Uncertain Variables:}
\begin{IEEEdescription}[\IEEEusemathlabelsep\IEEEsetlabelwidth{$Ptt$}]
  \item[$\tilde P_{d,t}^{\text{L/U}}$] Lower/Upper power bound of ADN $d$ at time slot $t$.
  \item[$\tilde E_{d,t}^{\text{L/U}}$] Lower/Upper energy bound of ADN $d$ at time slot $t$.
  \item[$\tilde{\xi}_{w,t}$] Forecasting error of renewable $w$.
  \item[$\tilde{\Omega}^{+/-}_{t}$] Total positive/negative renewable forecasting error at time slot $t$.
\end{IEEEdescription}
\subsection*{Decision Variables:}
\begin{IEEEdescription}[\IEEEusemathlabelsep\IEEEsetlabelwidth{$Ptt$}]
  \item[$P_{g,t}$] Active power of generator $g$ at time slot $t$.
  \item[$P_{g,s,t}$] Active power of the $s$-th segment of generator $g$ at time slot $t$.
  \item[$R_{g,t}^{\text{up/dn}}$]  Up-/Down-reserve of generator $g$ at time slot $t$. 
  \item[$\alpha_{g,t}^{+/-}$] Positive/Negative participation factor of generator $g$ at time slot $t$.
  \item[$P_{d,t}$] Active power of ADN $d$ at time slot $t$.
  \item[$R_{d,t}^{\text{up}}$] Up-reserve of ADN $d$ at time slot $t$.
  \item[$\alpha_{d,t}^{+}$] Positive participation factor of ADN $d$ at time slot $t$.
\end{IEEEdescription}
\subsection*{Notations:}
\begin{IEEEdescription}[\IEEEusemathlabelsep\IEEEsetlabelwidth{$Pttttxxt$}]
 \item[$\mathbf x, x$] Bold letters, e.g., $\mathbf x$, $\mathbf s$, $\mathbf z$, ${\boldsymbol\xi}$, and $\mathbf c$, denote vectors, and the corresponding non-bold letters denote their components.
 \item[$\tilde{\boldsymbol\xi},{\boldsymbol\xi}$] Letters with a tilde, e.g., $\tilde{\boldsymbol\xi}$, denote uncertain variables, and removing the tilde, e.g., ${\boldsymbol\xi}$, denotes their realizations.
 \item[${[}n{]}$] For any positive integer $n$, $[n]\triangleq \{1,2,...,n\}$.
 \item[$\| \cdot \|$, $\| \cdot \|_*$] Norm and dual norm of a vector.
 \item[$\mathbb{P}(\cdot)$] Probability of event $(\cdot)$.
 \item[$\mathbb{E}(\cdot)$] Expectation of uncertain variable $(\cdot)$.
 \item[$\mathbb I(\cdot)$] Indicator function, $\mathbb I(\cdot)=1$ if event $(\cdot)$ happens and $0$ otherwise.
\end{IEEEdescription}

\section{Introduction}\label{sect:intro}
\IEEEPARstart{R}{enewable} energies such as wind power and photovoltaics are increasingly integrated into the power system, the uncertain nature of which poses challenges to the system's operation. In this context, flexibilities of the demand-side distributed energy resources (DERs), e.g., electric vehicles (EVs), thermostatically controlled loads (TCLs), and distributed energy storage systems (DESSs), have received significant attention from grid operators worldwide, such as in Europe\cite{radeckeMarketsLocalFlexibility2019}, China\cite{huabeirule}, and the United States\cite{pjmrule}. Encouraging the numerous DERs to participate in the power system dispatch enhances the systems' security, economic efficiency, and ability to integrate renewables. Due to the large number but small capacities, DERs are organized as active distribution networks (ADNs) or aggregators and then treated as one equivalent dispatchable unit from the perspective of the transmission level. The dispatchable ranges of these DERs are also susceptible to uncertain factors, such as the arrival and departure time of EVs and the ambient temperature of TCLs. Therefore, modeling these uncertainties has become essential for decision-making in current and future power systems.

Uncertainty modeling has been widely discussed in the literature. Typical technical routes include scenario-based stochastic programming\cite{papavasiliou2013multiarea}, robust optimization\cite{ben2009robust}, and chance-constrained programming (CCP)\cite{charnes1963deterministic}. However, stochastic programming faces challenges in controlling risk levels. Robust optimization, on the other hand, tries to keep the strategy feasible under the worst realization of uncertain variables, which seems desirable for power system operation. Nonetheless, a major challenge of robust optimization lies in how to construct an appropriate uncertainty set \cite{bertsimasConstructingUncertaintySets2009}. Take the robust dispatch problem of EVs as an example for further explanation, where the EVs' arrival and departure times are treated as uncertain variables. It is common that, among all the uncertainty realizations, the latest arrival time is later than the earliest departure time. Consequently, if the uncertainty sets are constructed directly according to the maximum range of arrival and departure times, the robust dispatch problem ends up infeasible. Similar infeasibility phenomena are also observed in cases of other devices. In other words, it is sometimes necessary to allow a certain violation rate for constraints with uncertain variables to make the decision problem feasible. Hence, this paper focuses on CCP, which aims to do optimization while ensuring that the uncertain constraints are satisfied with a pre-set probability. Many studies have employed CCP to address uncertainties in power system decision-making \cite{8434318,huo2021chance,9814811}.

The general form of a chance constraint is:
\begin{equation}
  \mathbb{P}\left( {g_j}({\mathbf{x}},{\tilde{\boldsymbol{\xi }}}) \le 0,\forall j=1,2,...,m \right) \ge 1 - \epsilon,\label{eq:Gen_CC}
\end{equation}
where ${\mathbf{x}}\in \mathcal X$ is the vector of decision variables; ${\tilde{\boldsymbol{\xi }}}$ is the vector of uncertain variables with support set $\Xi$, and we denote its realizations by ${\boldsymbol{\xi }}$; the function $g_j\;(\forall j, 1\le j \le m)$ is assumed convex and lower semicontinuous in ${\mathbf{x}}$ for almost every ${\boldsymbol{\xi }}\in\Xi$; $\epsilon \geq 0$ is a pre-set risk level; and $m$ is the number of inequalities included in the chance constraint. If $m=1$, chance constraint \eqref{eq:Gen_CC} is an individual chance constraint (ICC), and otherwise, it is a joint chance constraint (JCC). 

Existing research has explored various approximations and reformulations to make ICCs or JCCs tractable. The sample-based approximation is a classic approach that approximates the chance constraints through the probability distribution on a certain number of samples \cite{luedtke2008sample}. Reference \cite{ruszczynski2002probabilistic} uses samples to approximately transform the CCP into a mixed integer programming problem, and reference \cite{jiaIterativeDecompositionJoint2021} approximates the CCP by iteratively checking combinations of constraint violation events based on the sample set. However, both the mixed integer programming and the computation of event combinations suffer from the curse of dimensionality. Besides, these methods have difficulties maintaining out-of-sample performance, limiting their application. Another popular approximation method is the Conditional Value at Risk (CVaR) approximation, a conservative convex approximation of the nonconvex chance constraint \cite{chenCVaRUncertaintySet2010}. The CVaR approximation transforms the chance constraint into an explicit deterministic constraint, thereby having a wide application in many fields \cite{liu2020CVaR,zhang2015distributed}. Nevertheless, the CVaR approximation faces the problem of over-conservativeness. In the extreme case, the CVaR approximation can even render an initially feasible chance constraint infeasible. Therefore, the CVaR approximation also has deficiencies.

Another factor that affects the approximation and reformulation of the CCP is whether the true distribution information $\mathbb P$ is known. In power system applications, we usually cannot access the true distribution $\mathbb P$ but only finite samples of the uncertain variables, which are subjected to different unknown but real data-generated probability distributions. Consequently, the decision-making should be robust to various probability distributions, i.e., should be distributionally robust \cite{mohajerin2018data,guo2018data}. In the context of this paper, distributional robustness arises in the form of distributionally robust chance constraints (DRCCs), which have been widely utilized in power system optimization problems to deal with the uncertain output of renewables \cite{ordoudis2021energy,yangTractableConvexApproximations2022,zhaiDistributionallyRobustJoint2022,dingDistributionallyRobustJoint2022}. The tractable approximation of DRCC, especially distributionally robust joint chance constraint (DRJCC), is a key focus in the related research. Many studies use the Bonferroni inequality to decompose the DRJCC into multiple distributionally robust ICCs and then reformulate them into tractable approximations \cite{ordoudis2021energy,yangTractableConvexApproximations2022,zhaiDistributionallyRobustJoint2022,dingDistributionallyRobustJoint2022}. However, it has been shown that even the best Bonferroni decomposition can be over-conservative \cite{zymlerDistributionallyRobustJoint2013}. Some literature also proposes to directly transform the DRJCCs, including an exact mixed-integer reformulation based on the Wasserstein metric \cite{xieDistributionallyRobustChance2021}, a linear reformulation based on the Wasserstein metric and CVaR approximation \cite{xieDistributionallyRobustChance2021}, and a linear matrix inequality reformulation based on moment information and CVaR approximation \cite{zymlerDistributionallyRobustJoint2013}. Still, mixed-integer reformulations and CVaR-based continuous reformulations have the same limitations as mentioned before.

This paper focuses on the application of a power system dispatch problem, including uncertainties from both renewable generations and dispatchable ranges of ADNs. Different resources in the system are affected by varying kinds and degrees of uncertainties, thus the probability of chance constraints should be prescribed individually. For instance, the reserve capacity of generators and line flows are affected to varying degrees by the uncertain forecast errors of renewable generations, while the flexibilities of ADNs are mainly affected by numerous uncertain end-user behaviors. Therefore, the system operator should set different risk levels for constraints associated with different resources, leading to the necessity of considering multiple chance constraints with different risk levels in the dispatch problem. In reference \cite{7332992}, different risk levels are considered for chance constraints on various generators and lines, while reference \cite{7828060} applies different risk levels to voltage constraints and line capacity constraints. For dispatchable units like ADNs, risk levels serve as an important metric for evaluating the uncertainty of their flexibilities. The system operator should evaluate the risk levels of all ADNs based on their historical deviations in executing dispatch signals and distinguish their uncertain flexibility with separate risk levels.

  We also notice that some studies employ one JCC that includes all constraints affected by uncertainty, e.g., reference \cite{pena2020dc}. This setting allows the system operator to evaluate the risk of constraint violations in the whole system using a single metric without distinguishing the different risk levels of various resources. In this case, if one resource is affected by a high level of uncertainty, i.e., has a high risk level, the system operator can not set a specific risk level for it but can only increase the risk level of the whole system to find a feasible dispatch strategy. Consequently, this dispatch strategy has a high risk level for all resources in the power system, leading to high deviations in the execution of dispatch signals and even compromising the safety of system operation. Hence, considering multiple chance constraints with varying risk levels is more suitable for dispatch problems.

The extension from a single chance constraint to multiple chance constraints is straightforward in traditional CVaR-based chance constraint reformulations. However, the aforementioned conservativeness and potential infeasibility concerns of the CVaR-based reformulation limit its application. Recently, the ALSO-X method, a less conservative chance constraint approximation method, was proposed \cite{jiangALSOXALSOXBetter2022}, addressing the limitations of CVaR. Nevertheless, the ALSO-X framework cannot be simply extended to multiple chance constraints because it inherently couples the JCC approximation with the solution of the optimization problem.

In order to solve the above shortcomings of the existing method and to obtain a more efficient and reliable dispatch strategy, this paper proposes a novel multiple JCC approximation method for uncertainty modeling in dispatch problems. The main contributions of this paper are threefold:
\begin{itemize}
  \item We propose a tractable approximation method for multiple JCCs in power system dispatch problems, which solves the over-conservativeness and potential infeasibility concerns of the conventional CVaR method. The proposed method is a significant extension of the ALSO-X approximation from dealing with a single JCC to multiple JCCs. This is particularly suitable for power system applications that require different risk levels for various resources.
  \item The proposed approximation method is then extended to multiple data-driven DRJCCs that fit the practical scenario of power system dispatch problems in which the distribution of uncertain variables is often inaccessible, but finite samples are easy to obtain.
  \item We propose a multi-period dispatch model for integrated transmission and distribution networks, including uncertainties from both renewable generations and flexibilities provided by ADNs. The impact of uncertainty on line flows, generators, and flexibility of ADNs over all time intervals is modeled as multiple DRJCCs, and the proposed approximation method is applied to convert it into a tractable formulation. The proposed method enables better feasibility and a wider adjustment range for cost and reliability than the conventional CVaR approximation.
\end{itemize}

The remainder of this paper is organized as follows. Section II introduces the proposed approximation method for multiple JCCs or DRJCCs. The multiperiod dispatch model is presented in Section III, using multiple DRJCCs to describe the uncertainties in renewable generations and ADN flexibilities. In Section IV, we first give two small examples to demonstrate the superiority and necessity of the proposed method over CVaR and the original ALSO-X and then implement numerical simulations for the multiperiod dispatch model on the IEEE test systems. Section V concludes this paper.


\section{Multiple JCCs Approximation Method}
This section first introduces the original ALSO-X approach for approximating a single JCC and then the proposed approximation method for multiple JCCs. We especially focus on deriving a specific form suitable for power system applications. Finally, we give a reformulation of the proposed approximation method for DRJCCs that match the practical conditions of power system decision-making.

\subsection{Original ALSO-X for Approximating a Single JCC and the Sample-based Counterpart for Power System Applications}
Let us consider the following CCP with a single JCC:
\begin{equation}
{f^*} = \min \left\{ {{{\mathbf{c}}^{\top}}{\mathbf{x}}:\mathbb{P} \left( {{g_j}({\mathbf{x}},{\tilde{\boldsymbol{\xi }}}) \le 0,\forall j \in [m]} \right) \ge 1 - \epsilon } \right\}, \label{eq:SJCCP_origin}
\end{equation}
where $\mathbf c$ is a constant vector. Even if the probability distribution of the uncertainty vector ${\tilde{\boldsymbol{\xi }}}$ is known, CCP \eqref{eq:SJCCP_origin} is generally intractable due to the probabilistic inequality that contains both decision and uncertain variables. To form a tractable approximation, CCP \eqref{eq:SJCCP_origin} is first equivalently transformed into a bilevel problem:
\begin{subequations}\label{eq:SJCCP_bilevel_eqv}
\begin{align}
  &\;\;\;\;\;\;\;\;\;{f^*} = \min f,\label{eq:SJCCP_bilevel_eqv:1}\\
  &\begin{array}{l}
    \text{s.t.}\;\;\left( {{{\mathbf{x}}^*},{s^*}( \cdot ),{z^*}( \cdot )} \right) = \\
    \;\;\;\;\;\;\;\arg \min \left\{ \begin{array}{l}
    \mathbb{E}[z({\tilde{\boldsymbol{\xi }}})s({\tilde{\boldsymbol{\xi }}})]:\\
    {\mathbf{x}} \in {{\mathcal X}},0\le z( \cdot ) \le 1,s( \cdot ) \ge 0,\\
    \mathbb{E}[z(\tilde{\boldsymbol{\xi }})] \ge 1 - \epsilon ,{{\mathbf{c}}^{\top}}{\mathbf{x}} \le f,\\
    {g_j}({\mathbf{x}},{\tilde{\boldsymbol{\xi }}}) \le s({\tilde{\boldsymbol{\xi }}}),\forall j \in [m]
    \end{array} \right\},
    \end{array} \label{eq:SJCCP_bilevel_eqv:2}\\
  &\;\;\;\;\;\;\;\;\;\mathbb{P} \left( s^*({\tilde{\boldsymbol{\xi }}})=0 \right) \ge 1 - \epsilon,\label{eq:SJCCP_bilevel_eqv:3}
\end{align}
\end{subequations}
where $s(\cdot):\Xi \mapsto \mathbb R$ and $z(\cdot):\Xi \mapsto \mathbb R$ are auxiliary functional variables (functions of ${\tilde{\boldsymbol{\xi }}}$) and $\mathcal X$ is the deterministic feasible region of the decision-variable $\mathbf x$. In the lower-level problem \eqref{eq:SJCCP_bilevel_eqv:2}, the upper bound of ${g_j}({\mathbf{x}},{\tilde{\boldsymbol{\xi }}})$ is relaxed by $s({\tilde{\boldsymbol{\xi }}})$, while variable $z({\tilde{\boldsymbol{\xi }}})$ represents the activation percentage of relaxation $s({\tilde{\boldsymbol{\xi }}})$, and the objective is to minimize the expectation of the activated relaxation. The equivalence between problems \eqref{eq:SJCCP_origin} and \eqref{eq:SJCCP_bilevel_eqv} is ensured by requiring the expected activation of relaxed constraints $ \mathbb{E}[z(\tilde{\boldsymbol{\xi }})]$ to be no less than $1 - \epsilon$ in the lower-level problem \eqref{eq:SJCCP_bilevel_eqv:2} and constraining the probability that relaxations $s^*({\tilde{\boldsymbol{\xi }}})$ being zero to be no less than $1 - \epsilon$ in constraint \eqref{eq:SJCCP_bilevel_eqv:3}. The detailed proof can be found in \cite{jiangALSOXALSOXBetter2022}.

Problem \eqref{eq:SJCCP_bilevel_eqv:2} does not include the probabilistic inequality that contains both decision and uncertain variables. If problem \eqref{eq:SJCCP_bilevel_eqv:2} is solved so that $\mathbf x^*$ is determined, then the probabilistic inequality \eqref{eq:SJCCP_bilevel_eqv:3} should be evaluated to make sure the solution $\mathbf x^*$ is feasible to the original CCP \eqref{eq:SJCCP_origin}. Therefore, there are still two barriers to overcome to make the equivalent problem \eqref{eq:SJCCP_bilevel_eqv} solvable: 1) how to evaluate the expectation in \eqref{eq:SJCCP_bilevel_eqv:2} and the probability in \eqref{eq:SJCCP_bilevel_eqv:3}, and 2) how to handle the bilinear item $z({\tilde{\boldsymbol{\xi }}})s({\tilde{\boldsymbol{\xi }}})$ in the objective of the lower-level problem \eqref{eq:SJCCP_bilevel_eqv:2}.

Regarding the first barrier, we consider the actual situation of power systems: the true distribution of the uncertainty vector ${\tilde{\boldsymbol{\xi }}}$ is usually unknown, and finite samples are available. Thus, by assuming that the uncertainty vector ${\tilde{\boldsymbol{\xi }}}$ is uniformly distributed over the finite sample set, which constitutes the empirical distribution, we reach the following sample-based counterpart of \eqref{eq:SJCCP_bilevel_eqv}:
\begin{subequations}\label{eq:SJCCP_bilevel_sample}
  \begin{align}
    &\;\;\;\;\;\;\;\;{f^*} = \min f,\label{eq:SJCCP_bilevel_sample:1}\\
    &\begin{array}{l}
      \text{s.t.}\;\left( {{{\mathbf{x}}^*},\mathbf{s^*},\mathbf{z^*}} \right) = \\
      \;\;\;\;\;\;\;\;\arg \min \left\{ \begin{array}{l}
          \frac{1}{n}\sum\nolimits_{i=1}^n {{z_i s_i}}:\\
          {\mathbf{x}} \in {{\mathcal X}},\mathbf 0\le{\mathbf{z}}\le \mathbf 1,{\mathbf{s}}\ge \mathbf 0,\\
          \frac{1}{n}\sum\nolimits_{i=1}^n {{z_i}}\ge 1 - \epsilon ,{{\mathbf{c}}^{\top}}{\mathbf{x}} \le f,\\
          {g_j}({\mathbf{x}},{\boldsymbol{\xi }_i}) \le s_i,\forall j \in [m],i\in[n]
        \end{array} \right\},
        \end{array} \label{eq:SJCCP_bilevel_sample:2}\\
    &\;\;\;\;\;\;\;\;\frac{1}{n}\sum\nolimits_{i=1}^n {\mathbb{I}(s_i^* = 0)}  \ge 1 - \epsilon,\label{eq:SJCCP_bilevel_sample:3}
  \end{align}
  \end{subequations}
  where $\boldsymbol{\xi }_i$ denotes the $i$-th sample of the uncertainty vector ${\tilde{\boldsymbol{\xi }}}$, and $n$ denotes the number of samples. The expectations and probability in problem \eqref{eq:SJCCP_bilevel_eqv} are converted into explicit linear expressions in problem \eqref{eq:SJCCP_bilevel_sample}, so that the first barrier is solved. 

However, there is still a bilinear term $z_i s_i$, rendering the lower-level problem \eqref{eq:SJCCP_bilevel_sample:2} non-convex. To obtain a tractable approximation, the original ALSO-X enforces ${\mathbf{z}} \equiv \mathbf 1$, i.e., assuming all the relaxations $s_i\;(\forall i\in[n])$ are fully activated, so that problem \eqref{eq:SJCCP_bilevel_sample} turns into:
\begin{subequations}\label{eq:SJCCP_bilevel_ALSO}
  \begin{align}
    &{f^*} = \min f, \text{  s.t.}\label{eq:SJCCP_bilevel_ALSO:1}\\
    &\left( {{{\mathbf{x}}^*},\mathbf{s^*}} \right) =\arg \min \left\{ \begin{array}{l}
        \frac{1}{n}\sum\nolimits_{i=1}^n {{s_i}}:\\
      {\mathbf{x}} \in {{\mathcal X}},{\mathbf{s}}\ge 0,{{\mathbf{c}}^{\top}}{\mathbf{x}} \le f\\
      {g_j}({\mathbf{x}},{\boldsymbol{\xi }_i}) \le s_i,\forall j \in [m],i\in[n]
      \end{array} \right\},\label{eq:SJCCP_bilevel_ALSO:2}\\
    &\frac{1}{n}\sum\nolimits_{i=1}^n {\mathbb{I}(s_i^* = 0)}  \ge 1 - \epsilon. \label{eq:SJCCP_bilevel_ALSO:3}
  \end{align}
  \end{subequations}

For any given $f$, the lower-level problem \eqref{eq:SJCCP_bilevel_ALSO:2} is convex and can be easily solved. Once \eqref{eq:SJCCP_bilevel_ALSO:2} is solved, it is straightforward to evaluate condition \eqref{eq:SJCCP_bilevel_ALSO:3}. Hence, a bisection of the objective $f$ can be employed to find a solution to problem \eqref{eq:SJCCP_bilevel_ALSO}, which results in the ALSO-X algorithm:
\begin{figure}[htbp]
  \vspace{-0.2in}
  \renewcommand{\algorithmicrequire}{\textbf{Input:}}
  \renewcommand{\algorithmicensure}{\textbf{Output:}}
  \removelatexerror
  \begin{algorithm}[H]
    \caption{Original ALSO-X Algorithm \cite{jiangALSOXALSOXBetter2022}}\label{alg:origin_ALSOX}
    \begin{algorithmic}[1]
      \REQUIRE Tolerence level $\delta_1>0$, upper and lower bounds $f^{\text{U}}$ and $f^{\text{L}}$ of the objective value
      \WHILE {$f^{\text{U}}-f^{\text{L}}>\delta_1$}
      \STATE Let $f = (f^{\text{U}}+f^{\text{L}})/2$, and solve problem \eqref{eq:SJCCP_bilevel_ALSO:2} to get $\mathbf x^*$ and $\mathbf s^*$
      \STATE Let $f^{\text{U}}=f$ if $\frac{1}{n}\sum\nolimits_{i=1}^n {\mathbb{I}(s_i^* = 0)}  \ge 1 - \epsilon$; otherwise, let $f^{\text{L}}=f$
      \ENDWHILE
      \ENSURE A solution $\mathbf x^*$ and its objective $f$
    \end{algorithmic}
  \end{algorithm}
  \vspace{-0.2in}
\end{figure}

The obtained solution $\mathbf x^*$ is feasible to the sample-based JCC \eqref{eq:SJCCP_bilevel_ALSO:3}, and the corresponding objective $f$ is no less than the optimal objective of problem \eqref{eq:SJCCP_bilevel_sample} because ${\mathbf{z}}$ is forced to equal to $\mathbf 1$. Hence, the obtained solution $\mathbf x^*$ is a conservative approximation of the optimal solution to problem \eqref{eq:SJCCP_bilevel_sample}, but it has been proved to be less conservative than the solution of the CVaR approximation \cite{jiangALSOXALSOXBetter2022}.

The initial upper and lower bound values $f^{\text{U}}$ and $f^{\text{L}}$ are straightforward to determine in practical problems based on their physical meaning. If the approximation problem \eqref{eq:SJCCP_bilevel_ALSO} is infeasible, the objective function series obtained in the iteration will converge to the initial upper bound $f^{\text{U}}$.

\subsection{Multiple JCCs Approximation}\label{sect:ALSO-X MJCC}
The original ALSO-X algorithm can yield satisfactory results in many practical problems with a single JCC. Nevertheless, it is necessary to separately address the risk levels of different chance constraints associated with different resources by considering multiple JCCs in practical applications because they may be affected by varying degrees of uncertainty or have different importance to the system. Consider the following CCP with multiple JCCs:
\begin{subequations}\label{eq:MJCCP}
  \begin{align}
    &{f^*} = \min {{\mathbf{c}}^{\top}}{\mathbf{x}},\label{eq:MJCCP:1}\\
    \text{s.t.}\;\; &{\mathbf{x}} \in {{\mathcal X}},\label{eq:MJCCP:2}\\
    &\mathbb{P} \left( {{g_{l,j}}({\mathbf{x}},\tilde{\boldsymbol{\xi}}_l) \le 0,\forall j \in [m_l]} \right) \ge 1 - {\varepsilon _l},\forall l \in [M]. \label{eq:MJCCP:3}
  \end{align}
  \end{subequations}
where $M$ is the number of JCCs; $\tilde{\boldsymbol{\xi}}_l$ is the vector of uncertain variables in the $l$-th JCC; $m_l$ and $\epsilon_l\ge0$ are the number of individual inequalities and the chosen risk level of the $l$-th JCC, respectively. 

Intuitively, the original ALSO-X algorithm can be extended to multiple JCCs by changing the judgment on whether \eqref{eq:SJCCP_bilevel_ALSO:3} holds, indicating that JCC \eqref{eq:SJCCP_bilevel_eqv:3} holds, to the judgment on whether all JCCs in \eqref{eq:MJCCP:3} hold simultaneously. This intuitive extension, however, can be very conservative compared to the original ALSO-X for a single JCC, particularly when each JCC is assigned a distinct confidence level. This limitation arises from the inability to distinguish the confidence levels of different JCCs in the lower-level problem \eqref{eq:SJCCP_bilevel_ALSO:2}.

Recall that the original ALSO-X approximation falls into this dilemma because it imposes $\mathbf z \equiv  \mathbf 1$ so that the constraint accounting for the risk level $\epsilon$ is removed from the sample-based equivalent form \eqref{eq:SJCCP_bilevel_sample:2}. Therefore, we still consider the sample-based counterpart of CCP \eqref{eq:MJCCP} with multiple JCCs as follows:
\begin{subequations}\label{eq:MJCCP_bilevel_sample}
  \begin{align}
    &{f^*} = \min f, \text{  s.t.}\label{eq:MJCCP_bilevel_sample:1}\\
    &\begin{array}{l}
      \left( {{{\mathbf{x}}^*},\mathbf{s^*},\mathbf{z^*}} \right) = \\
        \arg \min \left\{ \begin{array}{l}
          \frac{1}{M}\sum\nolimits_{l=1}^M\frac{1}{n_l}\sum\nolimits_{i=1}^{n_l} {{z_{l,i} s_{l,i}}}:\\
          {\mathbf{x}} \in {{\mathcal X}},\mathbf 0\le{\mathbf{z}}\le \mathbf 1,{\mathbf{s}}\ge \mathbf 0,{{\mathbf{c}}^{\top}}{\mathbf{x}} \le f,\\
          \forall l\in [M]:\\
          \frac{1}{n_l}\sum\nolimits_{i=1}^{n_l} {{z_{l,i}}}\ge 1 - \epsilon_l,\\
          {g_{l,j}}({\mathbf{x}},{\boldsymbol{\xi }_{l,i}}) \le s_{l,i},\forall j \in [m_l],i\in[n_l]
        \end{array} \right\},
        \end{array} \label{eq:MJCCP_bilevel_sample:2}\\
    &\frac{1}{n_l}\sum\nolimits_{i=1}^{n_l} {\mathbb{I}(s_{l,i}^* = 0)}  \ge 1 - \epsilon_l,\forall l\in [M],\label{eq:MJCCP_bilevel_sample:3}
  \end{align}
\end{subequations}
where $n_l$ denote the number of samples of the $l$-th vector $\tilde{\boldsymbol{\xi}}_l$. Similar to the relationship between problems \eqref{eq:SJCCP_origin} and \eqref{eq:SJCCP_bilevel_sample}, problem \eqref{eq:MJCCP_bilevel_sample} also serves as the sample-based equivalent form of CCP \eqref{eq:MJCCP} by requiring that the rate of activated relaxation of each JCC $l$ is no less than $1-\epsilon_l$ in the lower-level problem \eqref{eq:MJCCP_bilevel_sample:2} and that the probability that each JCC $l$ has zero relaxation over the sample set is no less than $1-\epsilon_l$ in the upper-level constraint \eqref{eq:MJCCP_bilevel_sample:3}. The theoretical proof is also similar to what proves the equivalence between problems \eqref{eq:SJCCP_origin} and \eqref{eq:SJCCP_bilevel_sample} in reference \cite{jiangALSOXALSOXBetter2022}, thereby omitted here for simplicity.

Rather than forcing $\mathbf z \equiv  \mathbf 1$, we propose to use an alternating optimization method, i.e., iteratively solving one of $\mathbf z $ and $\mathbf s$ by fixing the other, to solve the lower-level problem \eqref{eq:MJCCP_bilevel_sample:2}. Although a similar process has been applied to solve problem \eqref{eq:SJCCP_bilevel_sample} in \cite{jiangALSOXALSOXBetter2022} to improve the approximation accuracy for a single JCC, this improvement is insignificant because the original ALSO-X algorithm can adjust the solution based on the evaluation of condition \eqref{eq:SJCCP_bilevel_ALSO:3}. For multiple JCCs, however, simply evaluating the probabilities of all JCCs cannot provide a direction of how to adjust the intermediate solutions to match all these different confidence levels, while the alternating optimization can achieve this by accounting for the confidence levels in \eqref{eq:MJCCP_bilevel_sample:2} separately. Therefore, using the alternating optimization method to solve \eqref{eq:MJCCP_bilevel_sample:2} is more applicable to multiple JCCs and is expected to output a more accurate solution than the presented intuitive extension of the original ALSO-X.

Initialize $k=0$ and $\mathbf z^0=\mathbf 1$. At iteration $k+1$, we first fix $\mathbf z = \mathbf z^{k}$ in the lower-level problem \eqref{eq:MJCCP_bilevel_sample:2} and solve for $\mathbf s^*$. Since $\mathbf z$ has been fixed, the constraints on $\mathbf z$ can be ignored, and then problem \eqref{eq:MJCCP_bilevel_sample:2} becomes the following convex optimization:
\begin{equation}
  \begin{array}{l}
  \left( {{{\mathbf{x}}^{k+1}},\mathbf{s}^{k+1}} \right) = \\
    \arg \min \left\{ \begin{array}{l}
      \frac{1}{M}\sum\nolimits_{l=1}^M\frac{1}{n_l}\sum\nolimits_{i=1}^{n_l} {{z_{l,i}^k s_{l,i}}}:\\
      {\mathbf{x}} \in {{\mathcal X}},{\mathbf{s}}\ge \mathbf 0,{{\mathbf{c}}^{\top}}{\mathbf{x}} \le f,\\
      \forall l\in[M], j \in [m_l],i\in[n_l]: \\
      {g_{l,j}}({\mathbf{x}},{\boldsymbol{\xi }_{l,i}}) \le s_{l,i}
    \end{array} \right\}.
  \end{array}\label{eq:MJCCP_ALSOX_solve_sx}
\end{equation}
Then, we fix $\mathbf s = \mathbf s^{k+1}$ and solve for $\mathbf z^*$. Since $\mathbf s$ is constant, the constraints on $\mathbf s$ can be ignored, and problem \eqref{eq:MJCCP_bilevel_sample:2} becomes the following linear problem:
\begin{equation}
  \mathbf{z}^{k+1}_l=
  \arg \min \left\{ \begin{array}{l}
    \frac{1}{n_l}\sum\nolimits_{i=1}^{n_l} {{z_{l,i} s_{l,i}^{k+1}}}:\\
    \mathbf 0 \le \mathbf z_l \le \mathbf 1,\\
     \frac{1}{n_l}\sum\nolimits_{i=1}^{n_l} {{z_{l,i} }}\ge 1-\epsilon_l
  \end{array} \right\}, \forall l \in [M].
  \label{eq:MJCCP_ALSOX_solve_z}
\end{equation}
The $(\mathbf x^{k},\mathbf s^{k},\mathbf z^{k})$ obtained by this iteration are guaranteed to be feasible solutions of problem \eqref{eq:MJCCP_bilevel_sample:2} because constraints in problem \eqref{eq:MJCCP_bilevel_sample:2} are always satisfied. The iteration stops when the following criterion is met:
\begin{equation}
  \left| \frac{1}{M}{\sum\nolimits_{l = 1}^M {\frac{1}{{{n_l}}}\sum\nolimits_{i = 1}^{{n_l}} {\left( {z_{l,i}^{k + 1}s_{l,i}^{k + 1} - z_{l,i}^ks_{l,i}^k} \right)} } } \right| < {\delta _2},\label{eq:MJCCP_ALSOX_criteria}
\end{equation}
where $\delta _2>0$ is another pre-set tolerance parameter. 

The objective sequence $\{\frac{1}{M}\sum\nolimits_{l = 1}^M {\frac{1}{{{n_l}}}\sum\nolimits_{i = 1}^{{n_l}} {z_{l,i}^{k }s_{l,i}^{k}} } \}$ generated by the iteration is guaranteed to converge because it is monotonically non-increasing and has a lower bound of 0. The converged solution $\mathbf{x}^{k}$ is feasible for JCC \eqref{eq:SJCCP_bilevel_eqv:3} if $\frac{1}{M}\sum\nolimits_{l = 1}^M {\frac{1}{{{n_l}}}\sum\nolimits_{i = 1}^{{n_l}} {z_{l,i}^{k }s_{l,i}^{k}} } =0$; otherwise, it is infeasible.

The above discussion is summarized as the following approximation algorithm for multiple JCCs:
\begin{figure}[htbp]
  \vspace{-0.2in}
  \renewcommand{\algorithmicrequire}{\textbf{Input:}}
  \renewcommand{\algorithmicensure}{\textbf{Output:}}
  \removelatexerror
  \begin{algorithm}[H]
    \caption{Approximation Algorithm for Multiple JCCs}\label{alg:modified_ALSO}
    \begin{algorithmic}[1]
      \REQUIRE Tolerence level $\delta_1,\delta_2$, upper and lower bounds of the objective value $f^{\text{U}}$ and $f^{\text{L}}$
      \WHILE {$f^{\text{U}}-f^{\text{L}}>\delta_1$}
      \STATE Let $f = {(f^{\text{U}}+f^{\text{L}})}/{2}$, $\mathbf s^{0} = \boldsymbol{\infty} $, $\mathbf z^{0} = \mathbf 1$, $k = 0$
      \REPEAT 
      \STATE Solve problem \eqref{eq:MJCCP_ALSOX_solve_sx} to obtain $\mathbf s^{k+1}$
      \STATE Solve problem \eqref{eq:MJCCP_ALSOX_solve_z} for all $l\in[M]$ to obtain $\mathbf z^{k+1}$
      \STATE Let $\Delta = \left| \frac{1}{M}{\sum\nolimits_{l = 1}^M {\frac{1}{{{n_l}}}\sum\nolimits_{i = 1}^{{n_l}} {\left( {z_{l,i}^{k + 1}s_{l,i}^{k + 1} - z_{l,i}^ks_{l,i}^k} \right)} } } \right|$, $\Gamma = \frac{1}{M}\sum\nolimits_{l = 1}^M {\frac{1}{{{n_l}}}\sum\nolimits_{i = 1}^{{n_l}} {z_{l,i}^{k + 1}s_{l,i}^{k + 1}} } $, and $k=k+1$
      \UNTIL $\Gamma = 0$ or $\Delta<\delta_2$
      \STATE Let $f^{\text{U}} = f$ if $\Gamma =0$; otherwise, let $f^{\text{L}} = f$
      \ENDWHILE
      \ENSURE A solution $\mathbf x^{k}$ and its objective $f$
    \end{algorithmic}
  \end{algorithm}
  \vspace{-0.2in}
\end{figure}

The obtained solution $\mathbf x^{k}$ is feasible to the multiple sample-based JCCs \eqref{eq:MJCCP_bilevel_sample:3}, and the corresponding objective $f$ is no less than the optimal objective of problem \eqref{eq:MJCCP_bilevel_sample} because the alternating optimization does not guarantee to output the global optimal solution. Hence, the obtained solution $\mathbf x^{k}$ serves as a conservative approximation of the optimal solution for problem \eqref{eq:MJCCP_bilevel_sample}. Nevertheless, it is less conservative than the solution obtained through the intuitive extension of the original ALSO-X method and, certainly, less conservative than the solution derived from the CVaR approximation.

We here provide three Remarks with respect to the proposed algorithm:
\newtheorem{remark}{Remark}
\begin{remark}
  The parameter $\Gamma$ represents the objective of the lower-level problem \eqref{eq:MJCCP_bilevel_sample:2} obtained in the inner loop. Because $\Gamma = 0$ implies that the optimal solution of \eqref{eq:MJCCP_bilevel_sample:2} has been found, the stopping criterion $\Gamma = 0$ is added to the inner loop to reduce the number of iterations.
\end{remark}
\begin{remark}
  In the inner loop, when the given value of the objective $f$ is large, it is very likely to reach $\Gamma = 0$ after only one iteration. Thus, an initially large $f$ can be reduced rapidly. This is because all the constraints in \eqref{eq:MJCCP_bilevel_sample:3} can be easily satisfied when the objective $f$ is large enough, similar to what happens in the original ALSO-X algorithm.
\end{remark}
\begin{remark}
  The number of iterations required for the outer loop is $O(\log(f^{\text{U}}-f^{\text{L}}))$. Hence, it is advised to estimate the initial $f^{\text{U}}$ and $f^{\text{L}}$ as accurately as possible to reduce computation time. Fortunately, this is straightforward in real power systems. For example, the initial lower bound $f^{\text{L}}$ can be estimated by solving the deterministic version of \eqref{eq:MJCCP} where the uncertain variables are fixed to their sample average value, because the sample average approximation is almost always more optimistic than the chance constraints with a high confidence level. The initial upper bound $f^{\text{U}}$ can be estimated by solving problem \eqref{eq:MJCCP} via the CVaR approximation. If the CVaR approximation is infeasible, the initial upper bound $f^{\text{U}}$ can be estimated conservatively high based on the physics of the power system, which can then be reduced rapidly, as highlighted in Remark 2.
\end{remark}

\subsection{Data-Driven DRJCCs Approximation}\label{sect:ALSO-X DRJCC}
We have proposed the approximation method for multiple JCCs, which is a sample-based form assuming that the support set $\Xi$ only includes a finite number of samples of the uncertain variables. The sample-based form is a data-driven model that uses the sample data to represent the unknown distribution in real power systems for decision-making. Samples of these uncertain parameters can be obtained in various ways. For example, forecast samples of renewable generations can be obtained by various forecasting methods, and samples of ADNs' boundary parameters can be acquired from historical data since there is no well-recognized forecasting method for these parameters. However, the availability of the uncertainty samples in practical multiperiod dispatch problems is usually limited. This is because samples from recent similar days are more representative with less sampling error for the current dispatch decision. For example, the North China Power Grid uses historical operational data from the past ten days as references to evaluate the performance of various dispatchable resources \cite{huabeirule}. Hence, exploring out-of-sample information is imperative, raising the need for DRJCC, which is capable of trading off the out-of-sample performance with controllable conservativeness of dispatch decisions.

An optimization problem with multiple DRJCCs can be written as
\begin{equation}
  \begin{array}{c}
    {f^*} = \min {{\mathbf{c}}^{\top}}{\mathbf{x}},\text{  s.t.  }{\mathbf{x}} \in {{\mathcal X}},\\
    \inf\limits_{\mathbb P\in{\mathscr{P}_l}} \mathbb{P} \left( {{g_{l,j}}({\mathbf{x}},\tilde{\boldsymbol{\xi}}_l) \le 0,\forall j \in [m_l]} \right) \ge 1 - {\varepsilon _l},\forall l \in [M],
    \end{array}\label{eq:MDRJCCP}
\end{equation}
where $\mathscr{P}_l$ denotes the ambiguity set of the probability distribution in the $l$-th DRJCC. Each DRJCC in \eqref{eq:MDRJCCP} corresponds to the requirement that the inner JCC be satisfied even under the worst-case distribution within $\mathscr P_l$.

The construction of ambiguity set $\mathscr{P}_l$ is essential for a DRJCC. To incorporate the available sample data and to ensure a tractable reformulation, we define $\mathscr{P}_l$ as the \emph{Wasserstein ball} centered on the empirical distribution:
\begin{equation}
  \mathscr{P}_l = \left\{ {\mathbb P|\mathbb P( {\tilde{\boldsymbol{\xi}}_l  \in \Xi_l } ) = 1, W(\mathbb P,{{\hat{\mathbb P}}_{l,n_l}}) \le \rho_l } \right\},\label{eq:W_ball}
\end{equation}
where ${\hat{\mathbb P}}_{l,n_l}$ denotes the empirical distribution, i.e., the discrete uniform distribution on the $n_l$ samples of the uncertainty vector $\tilde{\boldsymbol{\xi}}_l$, $W(\cdot,\cdot)$ calculates the Wasserstein distance of two distributions, and $\rho_l \ge 0$ is a pre-set radius parameter for controlling conservativeness of the $l$-th DRJCC. When \(\rho_l = 0\), the DRJCC reduces to a sample-based JCC, and an increase in \(\rho_l\) expands the size of the ambiguity set $\mathscr{P}_l$, thereby resulting in a more conservative optimization result. The choice of the Wasserstein radius \(\rho_l\) depends on the decision-maker's degree of confidence in the samples' representativeness. The more the decision-maker trusts the samples' representativeness, the less conservativeness is needed in the decision-making, and the smaller \(\rho_l\) should be set. Adjusting \(\rho_l\) can also yield results for different levels of conservativeness, offering options for the decision-maker. In practice, the decision-maker can determine \(\rho_l\) in a data-driven manner, e.g., via cross validation [19]. In a word, the decision-maker can solve the optimization problem with DRJCCs under various \(\rho_l\)-s to identify the critical \(\rho_l\) that corresponds to the maximum profit under acceptable constraint violation rates.

Inspired by the original ALSO-X framework that can be extended to a tractable approximation of the data-driven DRJCC by modifying the constraint on $g_j$ in \eqref{eq:SJCCP_bilevel_ALSO:2}, we extend the proposed approximation algorithm to multiple DRJCCs by modifying the constraint on $g_{l,j}$ in \eqref{eq:MJCCP_ALSOX_solve_sx}. 
To this end, “solve problem \eqref{eq:MJCCP_ALSOX_solve_sx}" in line 4 of the proposed Algorithm \ref{alg:modified_ALSO} should be replaced by solving the following optimization problem:
\begin{equation}
  \begin{array}{l}
    \left( {{{\mathbf{x}}^{k+1}},\mathbf{s}^{k+1}} \right) = \\
      \arg \min \left\{ \begin{array}{l}
        \frac{1}{M}\sum\nolimits_{l=1}^M\frac{1}{n_l}\sum\nolimits_{i=1}^{n_l} {{z_{l,i}^k s_{l,i}}}:\\
        {\mathbf{x}} \in {{\mathcal X}},{\mathbf{s}}\ge \mathbf 0,{{\mathbf{c}}^{\top}}{\mathbf{x}} \le f,\\
        \forall l\in [M],j \in [m_l],i\in[n_l]:\\
        {\overline{g}_{l,j}}({\mathbf{x}},{\boldsymbol{\xi }_{l,i}}) \le s_{l,i}
      \end{array} \right\},
    \end{array} \label{eq:MDRJCCP_ALSOX_solve_sx_general}
\end{equation}
where the function ${\overline{g}_{l,j}}(\forall l\in[M],j\in[m_l])$ is defined as:
\begin{equation}
  {\overline{g}_{l,j}}({\mathbf{x}},{\boldsymbol{\xi }}) \triangleq {\max _{\boldsymbol{\zeta }}}\left\{ {{g_{l,j}}({\mathbf{x}},{\boldsymbol{\zeta }}),\text{s.t.,}\left\| {{\boldsymbol{\zeta }} - {\boldsymbol{\xi }}} \right\| \le \rho_l } \right\}\label{eq:g_up_general}.
\end{equation}

The derivation process for this reformulation is similar to that for extending the original ALSO-X for a single JCC to DRJCC \cite{jiangALSOXALSOXBetter2022}, thereby omitted. Comparing problem \eqref{eq:MJCCP_ALSOX_solve_sx} with \eqref{eq:MDRJCCP_ALSOX_solve_sx_general}, it can be observed that the extension from JCC to DRJCC is reflected in modifying the value of \({g_{l,j}}({\mathbf{x}},{\boldsymbol{\xi }})\) at \(\boldsymbol{\xi } = \boldsymbol{\xi }_{l,i}\) to the maximum value of \({g_{l,j}}({\mathbf{x}},{\boldsymbol{\xi }})\) within a neighborhood of \(\boldsymbol{\xi }_{l,i}\), leading to higher conservativeness. When the Wasserstein radius \(\rho_l = 0\), problem \eqref{eq:MDRJCCP_ALSOX_solve_sx_general} reduces to problem \eqref{eq:MJCCP_ALSOX_solve_sx}.

In power system dispatch problems, power flow constraints are usually linear \cite{zhaiDistributionallyRobustJoint2022,maghami2023two} and the flexibilities of ADNs can also be described by linear constraints \cite{silva2018estimating,kalantar2019characterizing}. Consequently, $g_{l,j} (\forall l\in[M],j\in[m_l])$ can be specified as a bi-affine form:
\begin{equation}
  {g_{l,j}}({\mathbf{x}},{\boldsymbol{\xi }}) =\boldsymbol{\xi}^{\top} \mathbf{a}_{l,j}(\mathbf x)+b_{l,j}(\mathbf x)\label{eq:g_biaffine},
\end{equation}
where $\mathbf{a}_{l,j}(\mathbf x)$ and $b_{l,j}(\mathbf x)$ are affine vector and scalar functions of $\mathbf x$, respectively. In this case, ${\overline{g}_{l,j}} (\forall l\in[M],j\in[m_l])$ has a closed-form expression:
\begin{equation}
  {\overline{g}_{l,j}}({\mathbf{x}},{\boldsymbol{\xi }}) =\rho_l {\left\| {{{\mathbf{a}}_{l,j}}(\mathbf x)} \right\|_*} + \boldsymbol{\xi}^{\top} \mathbf{a}_{l,j}(\mathbf x)+b_{l,j}(\mathbf x)\label{eq:g_up_biaffine}.
\end{equation}

By employing the $\ell_1$-norm in Equation \eqref{eq:g_up_general}, the dual norm in Equation \eqref{eq:g_up_biaffine} becomes the $\ell_\infty$-norm. As a result, problem \eqref{eq:MDRJCCP_ALSOX_solve_sx_general} can be specified as the following linear optimization:
\begin{equation}
  \begin{array}{l}
    \left( {{{\mathbf{x}}^{k+1}},\mathbf{s}^{k+1}} \right) = \\
      \arg \min \left\{ \begin{array}{l}
        \frac{1}{M}\sum\nolimits_{l=1}^M\frac{1}{n_l}\sum\nolimits_{i=1}^{n_l} {{z_{l,i}^k s_{l,i}}}:\\
        {\mathbf{x}} \in {{\mathcal X}},{\mathbf{s}}\ge \mathbf 0,{{\mathbf{c}}^{\top}}{\mathbf{x}} \le f,{\mathbf{v}_l}\ge \mathbf 0,\forall l\in [M],\\
        \|\mathbf{a}_{l,j}(\mathbf x)\|_\infty\le v_{l,j},\forall l\in [M],j \in [m_l],\\
        \forall l\in [M],j \in [m_l],i\in[n_l]:\\
        \rho_l v_{l,j}+\boldsymbol{\xi}_{l,i}^{\top} \mathbf{a}_{l,j}(\mathbf x)+b_{l,j}(\mathbf x) \le s_{l,i}
      \end{array} \right\}.
    \end{array} \label{eq:MDRJCCP_ALSOX_solve_sx}
\end{equation}
where $\|\cdot\|_\infty$ canculates the $\ell_\infty$-norm of a vector, and $\mathbf v_{l}\;(\forall l\in [M])$ is a vector of auxiliary variables.

When there are non-affinely adjustable resources in the system or the nonlinear power flow equations is considered, the function ${{g}_{l,j}}$ may not conform to the bi-affine form as in Equation \eqref{eq:g_biaffine}. In this case, the function ${\overline{g}_{l,j}}$ may not have a closed form expression like \eqref{eq:g_up_biaffine}, and it is necessary to derive ${\overline{g}_{l,j}}$ via Equation \eqref{eq:g_up_general} based on the specific ${{g}_{l,j}}$. Various ${\overline{g}_{l,j}}$ can be derived from different choices of ${{g}_{l,j}}$ (see \cite{ben2009robust} for detail). When a closed-form ${\overline{g}_{l,j}}$ is not available, a practical solution is to approximate the original ${{g}_{l,j}}$ using the ${{g}_{l,j}}$ families with a closed-form ${\overline{g}_{l,j}}$, e.g., the bi-affine family \eqref{eq:g_biaffine}. Errors introduced by this approximation can be compensated through real-time balancing services in the power system.

\section{Multiperiod Dispatch with Uncertainty}\label{sect:multi_period_dispatch}
The proposed approximation method for multiple JCCs or DRJCCs can be generally applied to the uncertainty modeling in many decision-making tasks of power systems, such as dispatch, planning, and control. Here, we focus on a multiperiod dispatch problem for the application of the proposed method. Power system dispatch determines the power reference points, participation factors, and reserve capacities of each dispatchable unit at an hourly or sub-hourly timescale, such as one hour, 15 min, or 5 min. Dispatch problems are formulated day-ahead, intraday, or both, depending on the requirements of the grid operators. 

The dispatch model in this paper includes renewable energies and conventional generators as the energy sources and conventional generators and ADNs as the dispatchable resources. The flexibility of an ADN is provided by the DERs connected to it. Since the ADN's capability to shift demand relies heavily on storage-like resources, such as DESSs, EVs, and TCLs, with strong temporal coupling, considering multiple time slots in the dispatch problem is necessary \cite{jabrRobustMultiPeriodOPF2015}. 

Uncertainties exist not only in renewable generation but also in the flexibility ranges of DERs, as mentioned in Section \ref{sect:intro}. These uncertainties affect the level of required reserves of generators, as well as the dispatchable ranges of ADNs and line flows. Different generators and ADNs may belong to different operators who can decide their risk aversion levels based on their exclusive operation preferences. The acceptable risk of violating line flow limits also depends on the importance of the respective line to the transmission network. Hence, it is necessary to distinguish between different risk levels using multiple DRJCCs.

An ADN acts as one dispatchable unit from the perspective of the transmission level. The dispatchable power range of an ADN results from aggregating the flexibility of all DERs connected to the distribution grid towards the transmission-distribution interface (root node of the distribution topology), considering distribution network constraints such as voltage limits and line flow capacities. In this paper, we adopt the power-energy boundary model \cite{wen2022aggregate,wen2023improved}, an approximate aggregation model with a clear physical meaning to describe the dispatchable ranges of ADNs in the multiperiod dispatch problem. The deterministic form of the power-energy boundary model is as follows:
\begin{subequations}\label{eq:peb}
  \begin{align}
    &P_t^{\text{L}} \le {P_t} \le P_t^{\text{U}},&\forall t \in [T],\label{eq:peb:p}\\
    &E_t^{\text{L}} \le \sum\nolimits_{\tau  = 1}^t {{P_\tau }\Delta T}  \le E_t^{\text{U}},&\forall t \in [T],\label{eq:peb:e}
  \end{align}
  \end{subequations}
where $\Delta T$ is the length of each time slot, and $T$ is the number of slots in the dispatch time horizon; $P_t$ is the power variable at time $t$; $P_t^{\text{L}}$, $P_t^{\text{U}}$, $E_t^{\text{L}}$, and $E_t^{\text{U}}$ are the power and energy boundary parameters. 

The power-energy boundary model, which includes intertemporal coupling, describes the dispatchable range of an ADN by constraining its power and accumulated energy consumption at each time slot. The boundary parameters $P_t^{\text{L}}$, $P_t^{\text{U}}$, $E_t^{\text{L}}$, and $E_t^{\text{U}}$ in the ADN's aggregated flexibility model are calculated based on the individual operational flexibility models of each DER within the ADN, taking into account the distribution network constraints, e.g., line capacities and voltage limits. Depending on the calculation method employed, these boundary parameters may render the power-energy boundary model either an outer or inner approximation of the exact aggregation model. Detailed calculation procedures of the boundary parameters in the outer and inner approximation models can be found in our previous works \cite{wen2022aggregate} and \cite{wen2023improved}, respectively. The boundary parameters are treated as uncertain variables in the following dispatch model because they are affected by the operational ranges of all DERs in the ADN.

The distributionally robust chance-constrained multiperiod dispatch problem is formulated as:
\begin{subequations}\label{eq:dispatch}
    \begin{equation}
      \begin{array}{l}
        \min \sum\limits_{t = 1}^T {\left[ {\sum\limits_{g \in {{\mathcal G}}} {\left( {{C_{g,0}}\Delta T + \sum\limits_{s = 1}^{{S_g}} {{C_{g,s}}{P_{g,s,t}}\Delta T} } \right)+} } \right.} \\
        \;\;\;\;\;\;\;\;\;\;\;\;\left. {\sum\limits_{g \in {{\mathcal G}}} {\left( {C_g^ + R_{g,t}^{{\text{up}}} + C_g^- R_{g,t}^{{\text{dn}}}} \right)}  + \sum\limits_{d \in {{\mathcal D}}} {C_d^+ R_{d,t}^{{\text{up}}}} } \right],
        \end{array}\label{eq:dispatch:obj}
    \end{equation}
s.t. 
    \begin{equation}
      0 \le {P_{g,s,t}} \le P_{g,s}^{\text{U} },\forall t \in [T],g\in\mathcal G,s\in[S_g],\label{eq:dispatch:1}
    \end{equation}    
    \begin{equation}
      \sum\limits_{s = 1}^{{S_g}} {{P_{g,s,t}}}  = {P_{g,t}},\forall t \in [T],\forall g \in {{\mathcal G}},\label{eq:dispatch:2}
    \end{equation}    
    \begin{equation}
      \sum\limits_{g \in {{\mathcal G}}} {{P_{g,t}}}  + \sum\limits_{w \in {{\mathcal W}}} {{P_{w,t}}}  = \sum\limits_{d \in {{\mathcal D}}} {{P_{d,t}}}  + \sum\limits_{i \in {{\mathcal I}}} {P_{i,t}^{{\text{fix}}}} ,\forall t \in [T],\label{eq:dispatch:3}
    \end{equation}    
    \begin{equation}
      \sum\limits_{g \in {{\mathcal G}}} {\alpha _{g,t}^ + }    = 1,\forall t \in [T],\label{eq:dispatch:4}
    \end{equation}    
    \begin{equation}
      \sum\limits_{g \in {{\mathcal G}}} {\alpha _{g,t}^ - }+ \sum\limits_{d \in {{\mathcal D}}} {\alpha _{d,t}^ + }  = 1,\forall t \in [T],\label{eq:dispatch:5}
    \end{equation}    
    \begin{equation}
      \alpha _{g,t}^ + \ge 0,\alpha _{g,t}^- \ge 0, \forall g \in \mathcal G, \alpha _{d,t}^ + \ge 0, \forall d \in \mathcal D, \forall t \in [T], \label{eq:dispatch:6}
    \end{equation}    
    \begin{equation}
      {P_{g,t}} + R_{g,t}^{{\text{up}}} \le P_g^{\text{max} },\forall t \in [T],g \in {{\mathcal G}},\label{eq:dispatch:7}
    \end{equation}    
    \begin{equation}
      {P_{g,t}} - R_{g,t}^{{\text{dn}}} \le P_g^{\text{min} },\forall t \in [T],g \in {{\mathcal G}},\label{eq:dispatch:8}
    \end{equation}  
    \begin{equation}
      r_g^{{\text{dn}}}\Delta T \le {P_{g,t + 1}} - {P_{g,t}} \le r_g^{{\text{up}}}\Delta T,\forall t \in [T - 1],g \in {{\mathcal G}},\label{eq:dispatch:9}
    \end{equation}  
    \begin{equation}
      \begin{array}{l}
        \mathop {\inf }\limits_{\mathbb{P} \in {\mathscr{P}_g}} \mathbb{P}\left\{ {R_{g,t}^{{\text{dn}}} \ge \alpha _{g,t}^ - \tilde{\Omega} _t^ + ,R_{g,t}^{{\text{up}}} \ge  - \alpha _{g,t}^ + \tilde{\Omega} _t^ - ,{\forall t \in [T]}} \right\}\\
        \;\;\;\;\;\;\;\;\;\;\;\;\;\;\;\;\;\;\;\;\;\;\;\;\;\;\;\;\;\;\;\;\;\;\;\;\;\;\;\;\;\;\;\;\;\;\;\;\;\;\ge 1 - {\varepsilon _g},\forall g \in {{\mathcal G}},\end{array}\label{eq:dispatch:10}
    \end{equation}  
    \begin{equation}
      \begin{array}{l}
        \mathop {\inf }\limits_{\mathbb{P} \in {\mathscr{P}_d}} \mathbb{P}{}\left\{ {{P_{d,t}} \ge \tilde P_{d,t}^{\text{L}},\sum\nolimits_{\tau  = 1}^t {{P_{d,t}}\Delta T}  \ge \tilde E_{d,t}^{\text{L}},} \right.\\
        \;\;\;\;\;\;\;\;{P_{d,t}} + R_{d,t}^{{\text{up}}} \le \tilde P_{d,t}^{\text{U}},\sum\nolimits_{\tau  = 1}^t {\left( {{P_{d,t}} + R_{d,t}^{{\text{up}}}} \right)\Delta T}  \le \tilde E_{d,t}^{\text{U}},\\
        \;\;\;\;\;\;\;\left. {R_{d,t}^{{\text{up}}} \ge \alpha _{d,t}^ + \tilde \Omega _t^ + ,\forall t \in [T]} \right\} \ge 1 - {\varepsilon _d},\forall d \in {{\mathcal D}},
        \end{array}\label{eq:dispatch:11}
    \end{equation}  
    \begin{equation}
      \begin{array}{l}
        \inf \limits_{\mathbb{P} \in {\mathscr{P}_l}} \mathbb{P}\left\{ { - P_l^{\text{U} } \le }{\sum\limits_{i \in {{\mathcal I}}} {{\psi _{l,i}}\left[ {\sum\limits_{g \in {{{\mathcal G}}_i}} {\left( {{P_{g,t}} - \alpha _{g,t}^ - \tilde \Omega _t^ +  - \alpha _{g,t}^ + \tilde \Omega _t^ - } \right)} } \right.} } \right. \\
        \;\;\;\;\;\;\;\;\;\;\;\;+ \sum\limits_{w \in {{{\mathcal W}}_i}} {\left( {{P_{w,t}} + {{\tilde \xi }_{w,t}}} \right)}  - \sum\limits_{d \in {{{\mathcal D}}_i}} {\left( {{P_{d,t}} + \alpha _{d,t}^ + \tilde \Omega _t^ + } \right)} \\
        \;\;\;\;\;\;\;\;\;\;\;\left. {\left. { - P_{i,t}^{{\text{fix}}}} \right]{ \le P_l^{\text{U} }},\forall t \in [T]} \right\} \ge 1 - {\varepsilon _l},\forall l \in {{\mathcal L}},
        \end{array}\label{eq:dispatch:12}
    \end{equation}  
    \end{subequations}
where the total positive and negative renewable forecasting errors, $\tilde{\Omega}^{+}_{t}$ and $\tilde{\Omega}^{-}_{t}$, are defined as:
\begin{equation}
  \tilde \Omega _t^ + \triangleq \max \left\{ {0,\sum\nolimits_{w \in {{\cal W}}} {{{\tilde \xi }_{w,t}}} } \right\},\tilde \Omega _t^ - \triangleq  \min \left\{ {0,\sum\nolimits_{w \in {{\cal W}}} {{{\tilde \xi }_{w,t}}} } \right\}.\nonumber
\end{equation}
Note that the directions of reserves are defined according to the perspective of the dispatchable unit: the up-reserve of a generator, $R_{g,t}^{\text{up}}$, means to increase generation, while the up-reserve of an ADN, $R_{d,t}^{\text{up}}$, means to increase the load.


  The most significant difference between the dispatch problem \eqref{eq:dispatch} and those in the literature lies in the asymmetrical modeling of participation factors and reserves. In references \cite{yangTractableConvexApproximations2022,jabrRobustMultiPeriodOPF2015,maghami2023two}, conventional generators are assumed to be the only reserve provider, so the up-/down-reserves and positive/negative participation factors are symmetrical. Since we have included ADNs as dispatchable units in the dispatch model, the reserves provided by ADNs should also be included. However, providing down-reserves requires the ADN's reference energy trajectory to be higher than its lower energy boundary, which leads the ADN to consume more energy. In typical scenarios where the energy cost is significantly higher than the reserve cost, if reserves and participation factors of ADNs are modeled symmetrically like what is usually done in the case of conventional generators, the optimal down-reserves of ADNs should be zero. Consequently, the optimal participation factors and up-reserves of ADNs are also zero. Hence, we propose the asymmetrical modeling of participation factors and reserves in \eqref{eq:dispatch:5} and \eqref{eq:dispatch:11}, assuming that ADNs only provide up-reserves so that the optimal up-reserve of ADNs can be nonzero to reduce the down-reserve provided by generators, which enables a more economically efficient dispatch result.

  The objective function \eqref{eq:dispatch:obj} includes the piecewise linearized energy costs of all generators and the reserve costs of all generators and ADNs over the entire optimization time horizon. There is no need for additional binary variables for the piecewise linearization because the cost functions of generators are monotonically increasing and convex. Constraints \eqref{eq:dispatch:1} and \eqref{eq:dispatch:2} are produced by the piecewise linearization of generators' energy costs. Equation \eqref{eq:dispatch:3} enforces the system-wide power balance. Constraints \eqref{eq:dispatch:4}-\eqref{eq:dispatch:6} ensure the adjustments in the power of generators and ADNs balance the renewable forecasting errors. Constraints \eqref{eq:dispatch:7}-\eqref{eq:dispatch:9} give the adjustable power range and ramping ability of each generator. DRJCCs \eqref{eq:dispatch:10}-\eqref{eq:dispatch:12} capture the probabilistic constraints on each generator, ADN, and transmission line, respectively. Constraint \eqref{eq:dispatch:10} enforces that the reserve activations requested from any generator $g$ over the optimization time horizon need to be within the corresponding reserve capacities $R_{g,t}^{\text{up}}$ and $R_{g,t}^{\text{dn}}$ with a probability $1-\epsilon_g$. Constraint \eqref{eq:dispatch:11} includes not only the reserve activation limits but also the dispatchable ranges of ADNs because the boundary parameters are uncertain variables. Finally, constraint \eqref{eq:dispatch:12} represents the probabilistic line flow limits. It is observed in constraints \eqref{eq:dispatch:10}-\eqref{eq:dispatch:12} that different resources in the system are affected by varying kinds and degrees of uncertainties. Hence, the system operator should set different risk levels for different resources. For example, the risk level of an ADN can be set based on its historical deviations from the dispatch strategy, assigning a high risk level for those ADNs with large deviations, and the risk levels of transmission lines can be determined according to their different importance in the system.

  The multiperiod dispatch problem \eqref{eq:dispatch} is of the form described in \eqref{eq:MDRJCCP}: the objective \eqref{eq:dispatch:obj} is linear, constraints \eqref{eq:dispatch:1}-\eqref{eq:dispatch:9} constitute a deterministic polytopic feasible region, viewed as the $\mathcal X$ in \eqref{eq:MDRJCCP}, and constraints \eqref{eq:dispatch:10}-\eqref{eq:dispatch:12} are the DRJCCs, where the inner inequalities are bi-affine with respect to the decision variables and uncertain variables. Therefore, the proposed approximation method in Subsections \ref{sect:ALSO-X MJCC} and \ref{sect:ALSO-X DRJCC} can be utilized to solve the multiperiod dispatch problem \eqref{eq:dispatch}. After transforming problem \eqref{eq:dispatch} into a tractable form using the proposed method, the scale of the problem is proportional to the scale of the transmission system and the number of uncertainty samples. Since the reformulated model is linear and, as noted at the beginning of Subsection \ref{sect:ALSO-X DRJCC}, does not require a significantly large dataset of samples, it will not lead to potential computational barriers.

  \section{Case Studies}
  In this section, we first test the effectiveness of the proposed JCCs approximation method using two small examples and then demonstrate its practicability by applying it to solve the proposed multiperiod dispatch model. All simulations are performed on a laptop with an 8-core Intel i7-1165G7 CPU and 32-GB RAM, programmed by \texttt{MATLAB}\cite{MATLAB}, and solved by \texttt{Gurobi}\cite{gurobi}.
  \subsection{Example 1: Upper and Lower Boundary Uncertainties}
  This example aims to compare the feasibility of ALSO-X with that of CVaR \cite{xieDistributionallyRobustChance2021} in the case of upper and lower boundary uncertainties. Although the advantages of the original ALSO-X over CVaR have been initially shown in \cite{jiangALSOXALSOXBetter2022}, the specific case of upper and lower boundary uncertainties, which represents the most simplified version of the ADN model in the dispatch problem, has yet to be discussed. We focus our discussion on JCC here rather than DRJCC to make the results clearer, because JCC is a simplified case of DRJCC, as stated in Subsection \ref{sect:ALSO-X DRJCC}. 
  
  Let us consider the following CCP with a single JCC:
  \begin{equation}
    \min f = x, \text{\;\;s.t.\;\;} \mathbb{P}(\tilde{\xi}^{\text{L}} \le x \le \tilde{\xi}^{\text{U}})\geq 1-\epsilon.\label{eq:toy1}
  \end{equation}
  Assume that we have five available samples of the uncertain variables $(\tilde{\xi}^{\text{L}},\tilde{\xi}^{\text{U}})$:
  \begin{equation}
    \left[ {\begin{array}{c}
      {\xi _{1,2,\cdots,5}^\text{L}}\\
      {\xi _{1,2,\cdots,5}^\text{U}}
      \end{array}} \right] = \left[ {\begin{array}{*{10}{c}}
      1&2&3&4&5\\
      3&4&5&6&7
      \end{array}} \right].\nonumber
  \end{equation}
The relationship ${\xi}^{\text{U}}>{\xi}^{\text{L}}$ holds for each pair of the given samples, and the maximum sample of ${\xi}^{\text{L}}$ is larger than the minimum sample of ${\xi}^{\text{U}}$. Under this sample set, the feasible region of \eqref{eq:toy1} is non-empty when $0.4\le\epsilon<1$.
  
  Using the ALSO-X approximation ($f^\text{U}=8$, $f^\text{L}=0$, $\delta_1 = 10^{-4}$) and the CVaR approximation \cite{xieDistributionallyRobustChance2021} to solve problem \eqref{eq:toy1} under different $\epsilon$-s, the results are shown in TABLE \ref{tab:toy1_results}. The results given by ALSO-X are the true optimal solutions, while CVaR cannot find a feasible solution independent of the value for $\epsilon$ within $[0,1)$ due to its over-conservativeness. 
  \begin{table}[htbp]
    \centering
    \captionsetup{justification=centering, labelsep=newline}
    \caption{\textsc{The Optimal Objectives of Example 1 under Different $\epsilon$ and Different Approximation Methods}}
    \setlength{\tabcolsep}{1mm}{
    \begin{tabular}{cccccc}
      \toprule
      $\epsilon$& 0 & 0.2 & 0.4 & 0.6 & 0.8\\
      \midrule
      ALSO-X  & infeasible  & infeasible   & 3  & 2 & 1\\
      \midrule
      CVaR   & infeasible  & infeasible   & infeasible  & infeasible & infeasible\\
      \bottomrule
      \end{tabular}%
    }
    \label{tab:toy1_results}%
  \end{table}

  This shortcoming of CVaR limits its application, particularly in modeling the uncertain dispatchable ranges of ADNs since both the upper and lower boundaries of ADNs are uncertain variables. Conversely, since ALSO-X can yield less conservative results for a single JCC, the proposed approximation method derived from ALSO-X is also expected to have a wider feasible range. This expectation will be further validated by a comparison between the proposed approximation method and CVaR in Subsection \ref{sect:simu_ed} in the simulations of the multiperiod dispatch problem.
      
 \subsection{Example 2: Multiple JCCs}
 This example aims to verify whether the proposed approximation method can better control the risk levels of different JCCs, so as to output a more accurate result than the intuitive extension of the original ALSO-X. The intuitive extension of the original ALSO-X to multiple JCCs, as briefly described at the beginning of Subsection \ref{sect:ALSO-X MJCC}, refers to changing the judgment on whether the single JCC is satisfied in line 3 of the original ALSO-X algorithm (Algorithm \ref{alg:origin_ALSOX}) to the judgment on whether all JCCs are satisfied simultaneously, if so, let $f^{\text{U} } = f$; and otherwise, $f^{\text{L}} =f$.

 Let us consider the following CCP with two JCCs:
 \begin{subequations}\label{eq:toy2}
  \begin{align}
    &\min f = y_1 +2y_2, \label{eq:toy2:obj}\\
    \text{s.t.}\;\;&x_1+x_2+x_3+x_4=y_1+y_2,\label{eq:toy2:1}\\
    &0\le y_1 \le 2, y_2\ge 0,\label{eq:toy2:2}\\
    &\mathbb{P}(x_1 \ge \tilde{\xi}_1,x_2 \ge \tilde{\xi}_2)\geq 1-\epsilon_1,\label{eq:toy2:3}\\
    &\mathbb{P}(x_3 \ge \tilde{\xi}_3,x_4 \ge \tilde{\xi}_4)\geq 1-\epsilon_2,\label{eq:toy2:4}
  \end{align}
  \end{subequations}
  where $x_1$, $x_2$, $x_3$, $x_4$, $y_1$, and $y_2$ are decision variables; $\tilde{\xi}_1$, $\tilde{\xi}_2$, $\tilde{\xi}_3$, and $\tilde{\xi}_4$ are uncertain variables uniformly distributed on $[0,1]$, with 20 samples drawn for each via random sampling. 
  
\begin{figure}[!t]
    \centering
    \subfigure[Objective function.]{\includegraphics[width=2.9in]{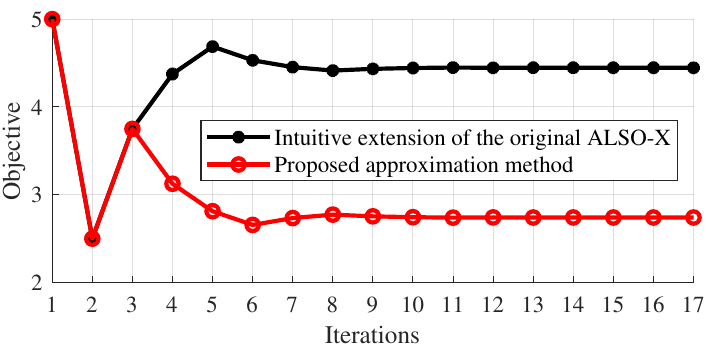}\label{fig:toy2_res:obj}}
    \subfigure[Violation rate of JCC \eqref{eq:toy2:3}.]{\includegraphics[width=3in]{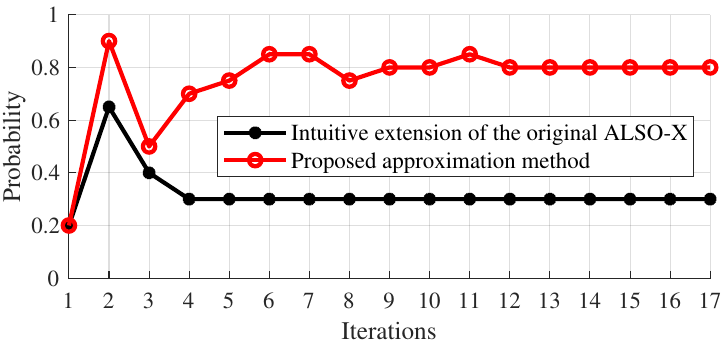}\label{fig:toy2_res:VP1}}
    \subfigure[Violation rate of JCC \eqref{eq:toy2:4}.]{\includegraphics[width=3in]{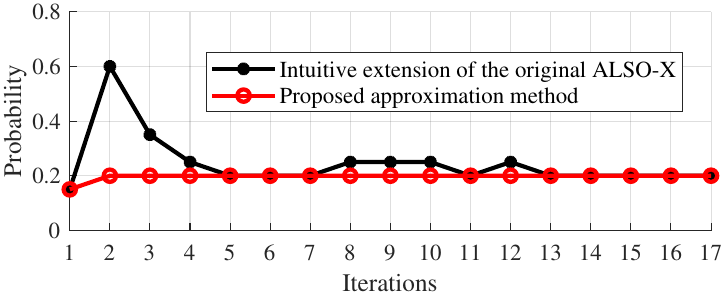}\label{fig:toy2_res:VP2}}
    \caption{Changes in the objective function and JCC violation rates during the iteration. The iteration numbers of the proposed approximation algorithm only count the outer loop.}
    \label{fig:toy2_res}
  \end{figure}

  We set $\epsilon_1=0.8$, $\epsilon_2=0.2$, $f^{\text{U}}=10$, $f^{\text{L}}=0$, and then use the intuitive extension of the original ALSO-X algorithm ($\delta_1 = 10^{-4}$) and the proposed approximation algorithm ($\delta_1 = \delta_2 = 10^{-4}$) to solve problem \eqref{eq:toy2}. Fig. \ref{fig:toy2_res} shows the changes in the objective function and the violation rates of JCC \eqref{eq:toy2:3} and \eqref{eq:toy2:4}, computed over the sample set, during the iteration. It is obvious that the proposed approximation method and the intuitive extension of the original ALSO-X converge to different results. In Fig. \ref{fig:toy2_res:obj}, the converged objective of the proposed method is significantly lower than that of the intuitive extension. The reason for this is indicated in Figs. \ref{fig:toy2_res:VP1} and \ref{fig:toy2_res:VP2}, where the proposed method precisely controls the violation rates of the two JCCs to the pre-set levels, i.e., $0.8$ and $0.2$, respectively, while the intuitive extension converges to an over-conservative result for JCC \eqref{eq:toy2:3}.

  Overall, the proposed approximation method outperforms the intuitive extension of the original ALSO-X in controlling the different risk levels of multiple JCCs. This advantage is necessary to take into account the risk aversion preferences of different stakeholders in power system decision-making.

  \subsection{Numerical Studies on the Multiperiod Dispatch Model}\label{sect:simu_ed}

  In this subsection, we consider a day-ahead 24-hour dispatch with an interval of $\Delta T = 1$ hour. The proposed approximation method is leveraged to reformulate and solve the multiperiod dispatch model with multiple DRJCCs. Additionally, we solve the problem using the CVaR approximated reformulation of DRJCCs \cite{xieDistributionallyRobustChance2021} in order to compare the solutions with respect to conservativeness and feasibility. 

  Simulations are carried out on an IEEE 30-node transmission network in the \texttt{Matpower} library. Here, we only include wind farms as renewable generation, but photovoltaics can be included in the same manner. Six wind farms are located at nodes 7, 14, 15, 26, 29, and 30, and six ADNs are located at nodes 4, 10, 12, 15, 17, and 19 in the transmission network. The reserve costs are $C_g^+ = C_g^- = \$2 /\text{MWh},\forall g\in \mathcal G$ and $C_d^+ = \$1 /\text{MWh},\forall d\in \mathcal D$. The installed capacity of each wind farm is 10 MW, and samples of output profiles comes from reference \cite{pinson2013wind}, which is forecasted based on realistic data in Denmark. Each ADN is a modified version of the IEEE 33-node distribution network, including 16 DER aggregators. In this study, we consider three types of DERs: EVs, heating, ventilating, and air conditioning systems (HVACs), and DESSs. In order to compare the conservativeness and feasibility characteristics under different scenarios, several cases are generated by adjusting the numbers and types of DERs assigned to each aggregator. However, since the conclusions of results in different cases are similar, we only select two typical cases to support our analysis, as follows:
\begin{itemize}
  \item \emph{Case 1}: 60 EVs, 40 HVACs, and 1 DESS for each aggregator;
  \item \emph{Case 2}: 30 EVs and 70 HVACs for each aggregator.
\end{itemize}

The operational parameters of EVs, HVACs, and DESSs are generated by random sampling based on their probability distributions, which are derived from numerous records of these devices in North China. Then, considering the LinDistFlow constraints of the distribution network \cite{baranOptimalSizingCapacitors1989}, the power-energy boundaries of ADNs are calculated using the \emph{flexibility optimization} algorithm (see \cite{wen2022aggregate} for the detailed calculation process). The calculated power and energy boundaries are input into the dispatch model as samples of the ADNs' uncertain variable of flexibility.

The number of uncertain samples is set to 20, representing 20 typical days. The risk level parameters $\epsilon_g$ and $\epsilon_d$ for DRJCCs \eqref{eq:dispatch:10} and \eqref{eq:dispatch:11} of each generator and ADN are set to 0.05, while the risk level $\epsilon_l$ for DRJCC \eqref{eq:dispatch:12} of each transmission line is set to 0.1. For simplicity, the Wasserstein radii of all DRJCCs are set identically, i.e., $\rho_g = \rho_l = \rho_d = \rho\;(\forall g\in \mathcal G,d\in \mathcal D,l\in \mathcal L)$, where $\rho$ can be adjusted to control the conservativeness. The initial lower bound $f^{\text{L}}$ of the objective in the ALSO-X algorithm is calculated by a deterministic dispatch based on the sample-average value of uncertain variables. The initial upper bound $f^{\text{U}}$ is calculated by the CVaR approximation if feasible, and otherwise, we set $f^{\text{U}} = 2f^{\text{L}}$.\footnote{In practical power system dispatch, the system's operational cost significantly surpasses the incremental costs generated by the conservative dispatch strategy against uncertainty for robustness. On the other hand, power systems' dispatch decisions usually have similar patterns over several consecutive days, resulting in relatively minor cost variations in daily dispatching strategies. Consequently, the initial upper bound $f^{\text{U}}$ can be approximated according to the deviation between the objective $f^*$ from the previous dispatch results and the lower bound $f^{\text{L}}$. For example, we can set the initial upper bound $f^{\text{U}}$ to $f^{\text{L}} + \alpha(f^* - f^{\text{L}})$, where $\alpha$ represents a coefficient that modulates the conservativeness of estimation. Typically, $\alpha$ can take values such as 2 or 3, thereby rendering the initial upper bound $f^{\text{U}}$ less conservative (lower) compared to a simple setting of $f^{\text{U}} = 2f^{\text{L}}$, so that the computational process can be accelerated. Nonetheless, this study adopts $f^{\text{U}} = 2f^{\text{L}}$ when the CVaR approximation is infeasible to demonstrate that the computational efficiency of the proposed method is sufficiently robust: even under highly conservative estimations of $f^{\text{U}}$, the computation time remains acceptable.} The convergence tolerances are $\delta_1 = 10^{-5}(f^{\text{U}}+f^{\text{L}})$ and $\delta_2 = 10^{-4}$.

In this simulation, we keep the input samples of the uncertain variables and risk level parameters unchanged and adjust the Wasserstein radius \(\rho\) of DRJCCs to compute the dispatch results. In the results for Case 1, CVaR and the proposed approximation methods yield feasible solutions in different ranges of $\rho$. We compute the changes in the total operation cost and \emph{out-of-sample reliability} of the two methods with the Wasserstein radius $\rho$. The out-of-sample reliability of a JCC is measured by counting the rate that all the inequalities that consitute the JCC are satisfied simultaneously under 100 test samples. These test samples are generated via the same random sampling process as how we generate the input sample set. The out-of-sample reliability of the DRJCCs depends on the specific DRJCC. Here, we only show the reliability results of DRJCC \eqref{eq:dispatch:11} for one ADN (the ADN at node 4) for simplicity.

\begin{figure}[!t]
  \centering
  \subfigure[Total operation cost.]{\includegraphics[width=3in]{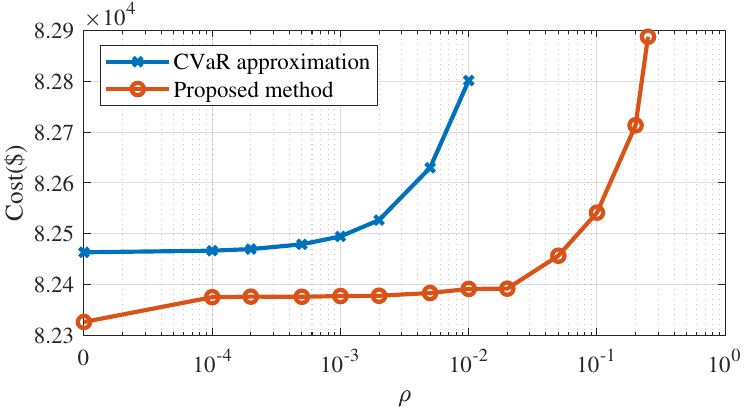}\label{fig:dispatch_case1_res:obj}}
  \subfigure[Out-of-sample reliability of the JCC for the ADN at node 4.]{\includegraphics[width=3in]{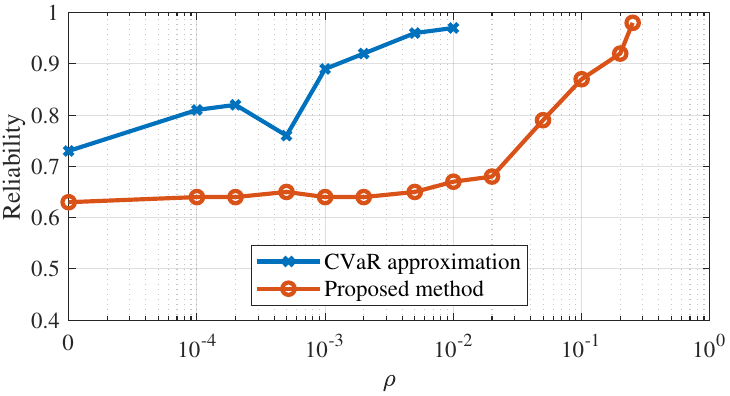}\label{fig:dispatch_case1_res:rlb}}
  \caption{Results of CVaR and the proposed approximation method for different Wasserstein radii $\rho$ in Case 1. The CVaR approximation leads to infeasibility when $\rho>0.01$.}
  \label{fig:dispatch_case1_res}
\end{figure}

Fig. \ref{fig:dispatch_case1_res} shows the cost and reliability of the two methods for different Wasserstein radii $\rho$. Both cost and reliability, in general, increase as $\rho$ increases with an exception when $\rho = 10^{-4}$. The CVaR approximation is always more conservative than the proposed approximation method at the same $\rho$. When $\rho$ increases, the CVaR approximation is feasible for $\rho \le 10^{-2}$, while the proposed approximation is feasible for $\rho \le  0.25 $. The proposed method ensures a wider adjustable range of $\rho$ than CVaR, thereby allowing a wider adjustable range of cost and reliability.
An interesting phenomenon can be observed if we closely examine Figs. \ref{fig:dispatch_case1_res:obj} and \ref{fig:dispatch_case1_res:rlb}. At certain points, such as with the proposed method at $\rho = 0.05$ and CVaR at $\rho = 0$, the proposed method can achieve higher reliability with a slightly lower cost. However, at other points, like with the proposed method at $\rho = 0.1$ and CVaR at $\rho = 0.001$, it is CVaR that can achieve higher reliability with a slightly lower cost. This variation is because different $\rho$ values actually lead to different settings of the ambiguity set, and the same sample set impacts the results of the two methods differently under these settings. In other words, the proposed method and CVaR are not consistently superior to each other in terms of cost and reliability under different settings of $\rho$. However, our main emphasis is that the proposed method has a broader adjustment range in cost and reliability than CVaR, indicating better feasibility. This advantage becomes much more essential in the subsequent Case 2, where CVaR is infeasible regardless of how the Wasserstein radius $\rho$ is adjusted, while the proposed method maintains a certain feasible range for $\rho$.

\begin{figure}[!t]
  \centering
  \subfigure[Accumulated energy plus up-reserve of the ADN and samples of the energy upper boundary.]{\includegraphics[width=3in]{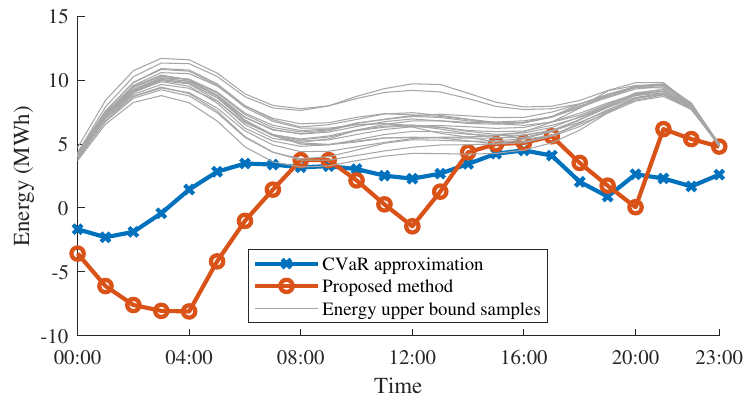}\label{fig:dispatch_case1_energ:upper}}
  \subfigure[Accumulated energy of the ADN and samples of the energy lower boundary.]{\includegraphics[width=3in]{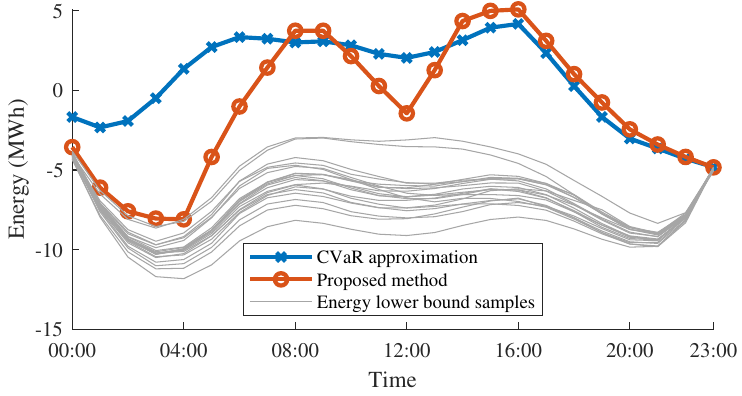}\label{fig:dispatch_case1_energ:lower}}
  \caption{Accumulated energy trajectories obtained by the two methods with the Wasserstein radius $\rho = 0.001$. The mean of all sample trajectories of the energy upper and lower bounds is subtracted from each displayed trajectory for clarity.}
  \label{fig:dispatch_case1_energy}
\end{figure}

Fig. \ref{fig:dispatch_case1_energy} further compares the dispatch results of the proposed method and CVaR for the ADN on node 4, where the average trajectory of all sample trajectories of the energy upper and lower bounds is subtracted from each displayed trajectory for clarity. Fig. \ref{fig:dispatch_case1_energ:upper} gives the accumulated energy trajectories for the ADN, including also the provision of up-reserve and the samples of the energy upper boundary, corresponding to the energy upper boundary constraint in \eqref{eq:dispatch:11}. Fig. \ref{fig:dispatch_case1_energ:lower} shows the accumulated energy trajectories of the ADN's power and the samples of the energy lower boundary, corresponding to the energy lower boundary constraint in \eqref{eq:dispatch:11}. Profiles of other constraints in \eqref{eq:dispatch:11} are omitted here for simplicity. The energy trajectory obtained by the proposed method is closer to the upper and lower boundary samples than that obtained by the CVaR method. At around 8:00 and 16:00, the energy trajectory obtained by the proposed method exceeds one upper boundary sample, which reflects the risk level of 0.05, i.e., 1/20. These results confirm that the proposed method is less conservative than CVaR and that it therefore better leverages the flexibility range of the ADN, resulting in a lower total cost of the proposed method compared to CVaR under the same Wasserstein radius $\rho$, as illustrated in Fig. \ref{fig:dispatch_case1_res:obj}.

\begin{figure}[!t]
  \centering
  \includegraphics[width=3.2in]{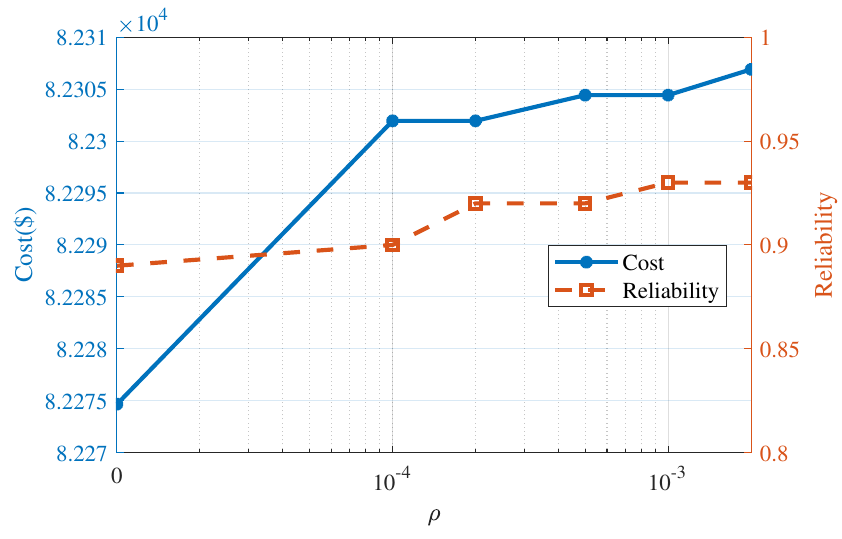}
  \caption{The cost and out-of-sample reliability of the proposed approximation method under varying Wasserstein radius $\rho$ in Case 2.}
  \label{fig:case2_cost_rlb}
\end{figure}

In Case 2, we still keep the input samples of the uncertain variables and risk level parameters unchanged and adjust the Wasserstein radius \(\rho\). The CVaR approximation can not yield a feasible solution no matter how we adjust the Wasserstein radius $\rho$, while the proposed approximation method stays feasible when $0 \le \rho \le 0.005$. Fig. \ref{fig:case2_cost_rlb} shows the cost and reliability of the proposed approximation method under varying Wasserstein radii $\rho$ in Case 2. It is observed that the adjustable range of cost and reliability in Case 2, as calculated by the proposed method, is narrower compared to Case 1. This is because the DER portfolio in Case 2 leads to more overlaps among the samples of the upper and lower boundaries in the ADN flexibility model. Nonetheless, the proposed method can, at least, maintain a certain range of feasibility, in contrast to the CVaR approximation, which remains infeasible regardless of adjustments in the Wasserstein radius $\rho$. This advantage in feasibility is important for the power system operator to avoid turning an originally feasible dispatch problem into an infeasible one.

To summarize, the proposed approximation method allows a more flexible tuning range for cost and reliability than the CVaR approximation. The proposed approximation method may find feasible solutions when CVaR is infeasible. The results demonstrate that the proposed method outperforms the conventional CVaR in the distributionally robust chance-constrained multiperiod dispatch problem.
\subsection{Scalability Test}
We further test the computation time of the proposed method in different systems. In addition to the study on the IEEE 30-node transmission system (with 6 ADNs), we perform the numerical experiments on an IEEE 57-node transmission system (with 13 ADNs) and on an IEEE 118-node transmission system (with 20 ADNs). The computation time on these three systems is 58.09 s, 248.64 s, and 998.63 s, respectively. Despite the observed increase in computation time with the system's scale, it is acceptable for the day-ahead 24-hour dispatch problem.

We also test the computation time in the case of intraday dispatch. The intraday dispatch problem considers much fewer time slots than the day-ahead dispatch, e.g., the North China Power Grid only considers the next 1 hour with a resolution of $\Delta T = 15$ minutes in the intraday dispatch, so $T=4$ in total. Under these settings of $T$ and $\Delta T$, the computation time of the dispatch model on the above three systems is 3.19 s, 10.84 s, and 39.82 s, respectively, which is still acceptable for intraday dispatch.


  \section{Conclusions}

  This paper proposes an approximation method to solve problems with multiple JCCs or DRJCCs. The proposed method is a significant modification of a recently proposed CCP approximation method called ALSO-X. The original ALSO-X is designed for a single JCC and cannot be intuitively extended to multiple JCCs because it couples the approximation of the JCC and the solution of the optimization problem. The proposed approximation method extends ALSO-X to multiple JCCs, which is necessary for power system decision-making because of the need to control risk levels of different uncertain resources separately.

  As an application scenario, we formulate a multiperiod dispatch model where the operational range of each generator, transmission line, and ADN is modeled as a DRJCC and thus can be solved with the proposed approximation algorithm. We highlight the discussion of the asymmetrical modeling of participation factors and reserves due to the properties of ADNs. Case studies show that the proposed approximation method can solve the over-conservativeness and potential infeasibility concerns of CVaR and is more precise in controlling the risk levels of different JCCs than the intuitive extension of the original ALSO-X method.

  A potential limitation of the proposed method is its reliance on samples. While benefiting from the model's efficiency, sample-based methods also suffer from the heavy dependence on whether the provided samples are representative enough for the real but unknown distributions. Therefore, developing forecasting methods for the uncertain variables in power systems, specifically forecasting the flexibility parameters of demand-side resources, becomes an important aspect of future research. Furthermore, the proposed approximation method can also be applied to a wide range of applications in power systems, including bidding, dispatch, planning, and control. Future works can develop decision-making tools for these specific problems to explore potential extension and generalization of the proposed method based on the distinct characteristics and requirements in these scenarios.


\small
\bibliographystyle{IEEEtran}

\bibliography{ref}


\end{document}